\documentclass[11pt]{article}

\usepackage{amsmath}
\usepackage{amsfonts}
\usepackage{enumerate}
\numberwithin{equation}{section}

\newtheorem{lemma}{Lemma}
\newtheorem{theorem}{Theorem}

\begin{document}

\author{Rafik Aramyan}

\title{RECONSTRUCTION OF CROFTON MEASURES FROM PROJECTIVE FINSLER METRICS
IN $\mathbf R^3$}

\maketitle

affiliation: {Institute of Mathematics Armenian Academy of Sciences
}

e-mail: {rafikaramyan@yahoo.com}

\abstract{A solution of Hilbert's fourth problem lead to integral equation of the type \emph{generalized cosine transform}. The present paper considers the solution  that integral equation by integral
geometry methods and propose an inversion formula for reconstruction of Crofton measures from projective smooth Finsler metrics in $\mathbf R^3$.}

\noindent subclass: {53C65,31A10}

\noindent  keywords: {integral geometry, integral equation, Finsler metrics, flag density}

\section{INTRODUCTION}

A solution of Hilbert's fourth problem interpreted as \emph{to construct and study all Finsler metrics in $\mathbf R^3$ whose
geodesics are straight lines} lead to integral equation of the type \emph{generalized cosine transform} (see \cite{Pog},\cite{Sch05}). It is the same to ask about construction of all projective (continuous and linearly additive) Finsler metrics in $\mathbf R^3$ (see \cite{Sch05},\cite{Alv03}).
The modern approaches make it clear that the
problem is at the basis of integral geometry, inverse problems and Finsler geometry.
First comes the A.V.Pogorelov-R.V.Ambartzumian-R.Aleksander result
(\cite{Amb76},\cite{Al},\cite{Pog}) stating that in two dimensional case every projective metrics can be obtained by Busemann's construction
with a measure in the space of lines in the plane. This essentially solves Hilbert's fourth problem in two dimensional case. In the smooth
case there is another approach via symplectic geometry due to Alvarez \cite{Alv}. Considerations on the construction of projective Finsler
metrics in two dimensional case are found by Ambartzumian and Oganyan in \cite{Amb98}.
One of the results in three dimensional case says that smooth Minkowski metrics in $\mathbf R^3$ are generated, via the cosine transform by signed measures in the space of directions in $\mathbf R^3$. By removing the
translation invariance restriction Pogorelov has shown in \cite{Pog} that any smooth projective Finsler metrics in $\mathbf R^3$ is generated by signed measure in the space of planes.

We denote $\mathbf E$ - the space of planes in $\mathbf R^3$, $\mathbf
S^2$ - the unit sphere in $\mathbf R^3$ (the space of unit vectors), $\mathbf S_\omega\subset\mathbf S^2$ - the oriented great circle with pole at
$\omega\in\mathbf S^2$. We consider
locally finite signed measures $\mu$ in the space $\mathbf E$, which
posses densities with respect to the standard Euclidean motion
invariant measure, i. e. (see \cite{Amb87})
\begin{equation}\label{1.1}
\mu(de)=h(e)de=h(e)\,dp\,d\xi\end{equation}
where $de$ is an element of the
standard measure, $(p,\xi)$ is the usual parametrization of a plane $e$: $p$ is
the distance of $e$ from the origin $O$; $\xi\in \mathbf
S^2$ is the direction normal to $e$.
By $[x]$ we denote the bundle of planes containing the
point $x\in \mathbf R^3$. To define function $h_{x}$ on ${\bold S}^2$ we consider the restriction of $h$ onto $[x]$ as a function on the hemisphere, since  a direction  completely
determines a plane from $[x]$. Then we extend the restriction to ${\bold S}^2$ by symmetry. Thus $h_{x}(\xi)=h(e_\xi)$ for $\xi\in{\bold S}^2$, where $e_\xi\in [x]$ is the plane with normal $\xi$. Below $h_{x}$ we call the restriction of $h$ onto $[x]$.

\noindent In \cite{Pog}, Pogorelov showed the following result.

\begin{theorem}
If $H$ is a smooth projective Finsler metric in $\mathbf R^3$ ($H:\mathbf R^3\times{\bold S}^2\rightarrow[0,\infty)$), then there exists a uniquely determined locally finite signed measure $\mu$ in the space $\mathbf E$, with
continuous density function $h$, such that

\begin{equation}\label{1.2}
H(x,\Omega)=
\int_{\mathbf
S^2}|(\Omega,\xi)|\,h_{x}(\xi)\,d\xi \,\,\,\,\emph{for}\,\,\,\,(x,\Omega)\in \mathbf R^3\times{\bold S}^2.\end{equation}
Here $h_{x}$ is the restriction of $h$ onto $[x]$, $d\xi$ denotes the spherical Lebesgue measure on $\mathbf
S^2$.\end{theorem}
The measure $\mu$ is also called a Crofton measure for the Finsler metric $H$ (see \cite{Sch05},\cite{Alv03}).
Different from the approach in \cite{Pog}, Theorem 1 was proved by Schneider \cite{Sch67} using expansions
in spherical harmonics.
Thus for smooth projective Finsler metrics, Pogorelov's result establish the existence of a Crofton measure, in general not positive. It seems that the problem of inversion of \eqref{1.2} and construction of the measure $\mu$ by means of $H$ was not considered before (we refer to the survey \cite{Sch05}).

In this paper we propose an inversion formula for \eqref{1.2} and reconstruct $h$ by given $H$ using integral and stochastic geometry approach.

We are interested in the inversion of the transform
$
h \to H$
defined by \eqref{1.2} (i.e. we
want to recover $h$ from a given $H$).

Note that, in case
$\mu$ is a translation invariant (t.i.) measure on $\mathbf E$ (the measure $\mu$ can be decomposed: there exists a finite even measure $m$ on $\mathbf
S^2$ such that $d\mu=m(d\xi)\cdot dp=h(\xi)d\xi\cdot dp$), \eqref{1.2}
represents the  \emph{cosine transform} (the zonoid equation)
\begin{equation}\label{1.4}
H(\Omega)=
\int_{\mathbf
S^2}|(\Omega,\xi)|\,h(\xi)\,d\xi.\end{equation}
playing an important role in convexity (see
\cite{SW}). If
the even function $H$ is sufficiently smooth, then \eqref{1.4} has a unique continuous solution $h$ and the inversion formula going back to W. Blaschke is well known. A rigorous proof specifying differentiability
assumptions was given by R.Schneider.

Note that one cannot reduce the solution of \eqref{1.2}
to the solution of the zonoid equation for a fixed point, since in general that solutions are not consistent. The last statement follows from the following result. Let $x\in\mathbf R^3$ be a fixed point and $\mu$ be a measure in the space $\mathbf E$, with
continuous density function $h$, which is not translation invariant. Now we consider the translation invariant measure $\mu_1$ whose element has the form $d\mu_1=h_x(\xi)d\xi\cdot dp$. Let $H$ and $H_1$ are the images of $\mu$ and $\mu_1$ respectively by transform \eqref{1.2}. It is obvious that the restrictions of $H(x,\cdot)$ and $H_1(x,\cdot)$ onto the point $x$ coincide.

In the present paper the problem of finding the solution of  \eqref{1.2} we
reduce to find the solution of an other
integral equation. The concept of a flag density which was introduced
and systematically employed by R. V. Ambartzumian, in
\cite{Amb87},\cite{Amb90} will be of basic importance below. We repeat the
definition.

We consider the so-called {\it directed} flags (below just a flag).
A flag is a triad $f=(x,g,e)$, where $x$ is a point in $\mathbf R^3$
called the {\emph location} of $f$, $g$ is a directed line containing the
point $x$, and $e$ is an oriented plane, i.e. a plane with
specified positive normal direction, containing $g$. There are two
equivalent representations of a flag:
$$f=(x,\Omega,\Phi)\,\,\,\text{or}\,\,\, f=(x,\omega,\varphi),$$
where $\Omega\in\mathbf
S^2$ is the spatial direction of $g$,
$\Phi$ is the planar direction in $\mathbf S_\Omega$ define the rotation of $e$ around $g$, $\omega\in\mathbf
S^2$ is the normal of
$e$, and $\varphi$ is the planar direction in $\mathbf S_\omega$ coincides with the direction
of $g$.

We use locally finite signed measure $\mu(de)=h(e)de$ in the space $\mathbf E$ to define the following function in the
space of flags $\mathcal F$ (so-called \emph{sine-square transform})(see \cite{Amb94},\cite{Ara94})
\begin{equation}\label{1.5}
\rho(f)=
\frac12\int_{\mathbf
S^2}\sin^2\alpha(\xi,f)\,h_{x}(\xi)\,d\xi\,\,\,\,\,\,\emph{for}\,\,\,\,\,\,f\in\mathcal F .\end{equation}
Here $h_{x}$ is the restriction of $h$ onto $[x]$. To explain $\alpha(\xi,f)$ we write $f=(x,g,e)$. Then $\alpha(\xi,f)$ is the angle between $g$ and the trace $e_\xi\cap e$, where $e_\xi$ is the plane through the origin ortohogonal to $e$. If we represent $f=(x,\omega,\varphi)$ then
\begin{equation}\label{1.5.1}
\sin^2\alpha(\xi,f)=\cos^2(\varphi-\psi),\end{equation}
where $\psi$ is the
direction of the projection of $\xi$ into the plane of the flag $f$ and \eqref{1.5.1} does not
depend on the choice of the reference point on the plane of the flag
$f$.
The function $\rho$ defined on the space of flags $\mathcal F$
we call \emph{flag densities} of measure $\mu$.

If the flag density $\rho$ is sufficiently often differentiable, then \eqref{1.5} has a unique continuous solution $h$. In \cite{Ara94} by author of the present paper using integral and stochastic geometry
methods was found an inversion formula for \eqref{1.5} and reconstruct the density $h$ of sign measure $\mu$ in terms of its smooth flag density $\rho$.

The problem of finding the solutions of  \eqref{1.2} we
reduce to finding the flag density $\rho$ for which
\begin{equation}\label{1.6}
H(x,\Omega)=\frac1{\pi}
\int_{0}^{2\pi}\rho(x,\Omega,\Phi)\,d\Phi.\end{equation}
and using inversion formula of \eqref{1.5}.

\noindent Now we describe the inversion formula of \eqref{1.2}.
We need to give the definitions of certain partial derivatives of $\rho(f)$.
With each directed flag we associate three orthogonal axes through
$x$: by definition, axis $x_1$ coincides with the direction of $g$;
axis $x_2$ lies within $e$, is orthogonal to $x_1$ and
is directed into the right half of $e$ bounded by $
g$; axis $x_3$ is coincides with the positive normal to $e$. We
require that the axes $x_1$,$x_2$,$x_3$ form a left triad.

By $\rho'_\Phi(f)$ we denote the derivative of $\rho$ at $f$ which
corresponds to positive rotation of $f$ around the axis
${g}$. By definition, the positive rotation of the
space around the axis ${g}$ appears clockwise, when
we look in the direction of the axis ${g}$.

We will also denote by $\frac{\partial}{\partial_\eta P}$ the
partial derivative in the argument $x$ which is taken in the
direction $\eta\in \mathbf{S}^2$. In the special cases where $\eta$
coincides with directions of the axes $x_2$ or $x_3$ respectively, the
values of $\frac{\partial\rho}{\partial_\eta x}$ will be denoted as
$\rho'_y$ and $\rho'_n$ correspondingly.

Our main result is the following. Let $H$ be a smooth projective Finsler metric on $\mathbf R^3$. For every fixed $x\in\mathbf R^3$, we now
solve the zonoid equation for the function $H(x,\cdot)$. Since this function is of
class $C^\infty$, there exists a continuous even solution
$h(x,\cdot)$ on ${\bold S}^2$.

Then we put the function $h(x,\cdot)$ into equation \eqref{1.5} instead of $h_{x}(\cdot)$ and found the flag function $\rho$
\begin{equation}\label{1.7}
\rho(f)=\rho(x,\omega,\varphi)=
\frac12\int_{\mathbf
S^2}\cos^2(\varphi-\psi)\,h(x,\xi)\,d\xi\,\,\, \emph{for every}\,\,\,\, (\omega,\varphi),\end{equation}
which depends on $x\in\mathbf R^3$. The function $\rho$ is the flag density of the  Crofton measure for the projective Finsler metric $H$.

For a given plane $e$ and a point
$x\in e$, the corresponding so-called bundle of flags we denote by
$$<x,e>=\{f:\,\,f=(x,g,e)\}.$$
We consider the  function so-called bundle mass at $x\in\mathbf R^3$:
\begin{equation}\label{1.8}
M(x)=\frac1{\pi}\int_0^{2\pi}\rho(x,\omega,\varphi)d\varphi=
\frac12\int_{\mathbf
S^2}h(x,\xi)\,d\xi.\end{equation}
Note that the first integral in \eqref{1.8} does not depend on $\omega\in\mathbf
S^2$.

By $\rho(\cdot)$
we denote the restriction of $\rho$ onto $<x,e>$. The notation
$\rho(\varphi)$ is reasonable since  $\varphi$ completely determines
a flag from $<x,e>$.

\begin{theorem}
\label{t1} Let $H$ be a smooth projective Finsler metric in $\mathbf R^3$, that is defined by transform \eqref{1.2}
of a sign measure $\mu$ in $\mathbf E$ with continuous density $h$ with respect to the invariant measure.
For a given plane $e$ the following representation is valid
\begin{equation}\label{1.9}
h(e)=M(x)+\frac1{2\pi}\int_0^{2\pi}
\frac{\partial^2M(x)}{\partial^2_\varphi x}\,d\varphi-\frac2\pi
\int_0^{2\pi}(\rho''_{\phi\phi}(\varphi)+2\rho''_{\phi
y}(\varphi)+\rho''_{yy}(\varphi))d\varphi.
\end{equation}
where $x$ is a point on $e$, $\rho$ is defined by \eqref{1.7}, $\rho(\cdot)$ is the restriction of $\rho$ onto $<x,e>$.
\end{theorem}
In \S 3 we present the
expression for $\rho'_\Phi(\varphi)$ in terms of the derivatives of
the function with respect to the parameters involved.

\section {FLAG DENSITIES AND MEASURES IN THE SPACE OF PLANES}
 \label{sec1}

Equation \eqref{1.5} naturally emerges in
Integral Geometry in the theory of recovering of a
measure defined in the space of planes in $\mathbf R^3$ from the
known flag density (see \cite{Amb94},\cite{Ara94}).

\noindent It was proved in \cite{Ara94} that equation \eqref{1.5} has the unique solution in the class of continuous functions
and found an inversion formula. Here we present a short version of the proof of the formula for
completeness.

Note that, in \cite{Ara92} the same problem was considered for the case
$\mu$ is a translation invariant measure on $\mathbf E$ (the measure $\mu$ can be decomposed: there exists a finite even measure $m$ on $\mathbf
S^2$ such that $d\mu=m(d\xi)\cdot dp=h(\xi)d\xi\cdot dp$). In
\cite{Ara92} (see also \cite{Ara09})  was obtained integral expression for the value of $m$ for
a spherical domain bounded by a piecewise smooth curve and and an inversion
formula was found for a case $m$ has a density, which first was found in \cite{Amb87}.

To invert equation \eqref{1.5} we do the following. Let $\rho$ be a
smooth function on the space of flags that is defined by
\eqref{1.5} of a measure $\mu$ in $\mathbf E$ possessing continuous
density $h$. We consider the restriction of $h$ onto the set
of planes tangent to a spherical domain. Then, by integral geometry
methods we find the integral of the restriction over a spherical
disc in terms of $\rho$. Using this integrals we find an
inversion formula for $h$.

We need some results from integral geometry.
Let $\mu$ be a sign measure on $\mathbf E$,
possessing a density $h$ with respect the invariant measure, i.e.
$\mu(de)=h(e)de$. Given a subset $A\subset\mathbf R^3$, by $[A]$ we
denote the set of planes, that intersect $A$. Let $\mathbf{B}$ be a
convex body with a sufficiently smooth boundary
$\partial\mathbf{B}$. By $k_1,k_2$ we denote the principal normal
curvatures of $\partial\mathbf{B}$ at a point
$s\in\partial\mathbf{B}$, and by $f_i=(s,g_i,t)$, $i=1,2$ we denote
the flag, where $t$ is the plane tangent to $\partial\mathbf{B}$ at
the point $s$, $g_i$ is the $i$-th principal direction of curvature
at $s\in\partial\mathbf{B}$. In \cite{Ara94}(see also \cite{Ara10}) was found the following representation for $\mu([\mathbf{B}])$.

\begin{theorem}
\label{t2} Let $\mu$ be a sign measure on $\mathbf E$, possessing a
density $h$ with respect the invariant measure. For any sufficiently smooth convex body $\mathbf{B}$
we have the following representation:
\begin{equation}\label{2.1}
\mu([\mathbf{B}])=(2\pi)^{-1}\int_{\partial\mathbf{B}}
[k_1\rho(f_2)+k_2\rho(f_1)]\,ds,\end{equation} where $\rho$ is the
flag density of $\mu$ defined by \eqref{1.5}, $ds$ is an area
element on $\partial\mathbf{B}$.
\end{theorem}

\noindent For the case
$\mu$ is a translation invariant measure on $\mathbf E$ the representation first was found in
\cite{Pa}.

\noindent We will need some further definitions. By
$\mathbf{S}(Q,R)$ we denote the sphere with center $Q$ and radius
$R$, $\mathbf{S}^2$ will stand for $\mathbf{S}(O,1)$.

\noindent Assume $h$ is the density of $\mu$. We define the following function on $\mathbf{S}(Q,R)$
$$h^*(s)=h(t(s)) \,\,\,\,\,\emph{for}\,\,\,\,\,\,s\in\mathbf{S}(Q,R),$$
where $t(s)$ is the plane tangent to $\mathbf{S}(Q,R)$ at $s\in\mathbf{S}(Q,R)$.

\noindent The measure $R^{-2}\, h^*(s)\,ds$, where $ds$ is an area
element on $\mathbf{S}(Q,R)$,we call the {\it conditional} measure
on $\mathbf{S}(Q,R)$, generated by $\mu$ (or $h$).

\noindent Let $\mathrm{A}$ be a geodezically convex domain on
$\mathbf{S}(Q,R)$, which is contained in some hemisphere. Let
$\mathbf{B}_\varepsilon\subset\mathbf{R}^3, \varepsilon>0$ be the
convex body bounded by \noindent $\mathbf{S}(Q,R+\varepsilon)$,
planes tangent to $\mathbf{S}(Q,R)$ at points of
$\partial\mathrm{A}$ and the conical surface with the vertex $Q$ and
base $\mathrm{A}$. In \cite{Ara94} was proved the following theorem.

\begin{theorem}
\label{t3} Let $\mu$ be a measure on $\mathbf{E}$ with continuous
density $h$ with respect the invariant measure, and let $\mathrm{A}\subset\mathbf{S}(Q,R)$ be an
open convex domain contained in some hemisphere. Then
\begin{equation}\label{2.2}
\frac1{R^2}\,\int_{\mathrm{A}}
h^*(s)\,ds=\lim_{\varepsilon\to0}\frac{\mu([\mathbf{B}_\varepsilon])-
\mu([\mathbf{B}_0])}\varepsilon.\end{equation}
\end{theorem}

\noindent Using Theorems 3 and 4, one can calculate the values of
the ``conditional measure" of various domains $A\subset\mathbf S^2$.

\noindent Now we calculate the conditional measure of
a spherical disc $\mathrm{A}\subset\mathbf{S}(Q,1)$ of spherical
radius $\alpha<\pi/2$ and using this result we
find an inversion formula for \eqref{1.5}. For translation invariant case the result was
obtained in \cite{Ara92} (see also \cite{Ara09}). Without loss of
generality, one can consider $Q$ as the origin $O$.
On $\mathbf S^2$ we consider usual spherical coordinates
$(\varphi,\nu)$. The center of the disc chosen for the pole.
We have:
$\partial A=\{l=(\varphi,\alpha):
0\leq\varphi<2\pi\}$.

\noindent A flag $f$ we call a tangent flag to $\mathrm{A}$ at
$l\in\partial A$, if $f$ is located at $l$, the plane of $f$ is
tangent to $\mathbf{S}^2$ at $l\in\partial\mathrm{A}$ and the positive normal of the plane of $f$ coincides with the outer
normal to $\mathbf S^2$ at $l\in\partial\mathrm{A}$, the line
of $f$ is tangent to $\partial \mathrm{A}$ at $l$ and the
direction on the line of $f$ corresponds to the motion along
$\partial\mathrm{A}$, which leaves (locally) $\mathrm{A}$ on the
left hand side. By
$<\partial\mathrm{A}>$ we denote the set of the tangent flags
of $\mathrm{A}$. The tangent flag at $l$ is unique and therefore
$\rho(l)$ is a reasonable notation for the value of $\rho$ at the tangent flag located at $l$. By $\rho(\cdot)$ we denote the restriction of a flag
function $\rho$ onto $<\partial\mathrm{A}>$.

\begin{theorem}
\label{t1} Let $\rho$ be a $C^2$ smooth function on the space of
flags in $\mathbf R^3$, that is defined by equation \eqref{1.5}
of a sign measure $\mu$ in $\mathbf E$ with density $h$ with respect to the invariant measure.
Then for any spherical disc $\mathrm{A}\subset\mathbf S^2$
\begin{multline}\label{2.3}
2\pi\int_{A} h^*(s)\,ds=\int_{\mathrm{A}}M(s)\,ds
+\int_{\mathrm{A}}\frac{\partial M(s)}{\partial n}\,ds
+\cos\alpha\int_0^{2\pi}M(l)\,d\varphi-\\
-\sin\alpha
\int_0^{2\pi}\rho'_\phi(l)\,d\varphi
-2\cos\alpha\int_0^{2\pi}\rho(l)\,d\varphi-
\sin\alpha\int_0^{2\pi}\rho'_y(l)\,d\varphi,\end{multline}
here $n$ is the outer normal direction to $\mathbf S^2$ at $s$, $d\varphi$ denotes the Lebesgue measure on $[0,2\pi)$.
\end{theorem}

{\emph{Proof of the Theorem 5}} According to Theorem 4 we have to calculate
$\mu([\mathbf{B}_\varepsilon])$ and $\mu([\mathbf{B}_0])$. In order
to apply Theorem 3, instead of $\mathbf{B}_\varepsilon$ and
$\mathbf{B}_0$ we consider their smooth versions
$\mathbf{B}_\varepsilon+\mathbf{S}(O,\delta)$ and
$\mathbf{B}_0+\mathbf{S}(O,\delta)$ for some $\delta>0$. We have
\begin{equation}\label{2.4}\mu([\mathbf{B}_\varepsilon])=
\lim_{\delta\to0}\mu([\mathbf{B}_\varepsilon+\mathbf{S}(O,\delta)]),
\,\,\,\mu([\mathbf{B}_0])=\lim_{\delta\to0}
\mu([\mathbf{B}_0+\mathbf{S}(O,\delta)]).\end{equation} We divide
the surfaces $\mathbf{B}_\varepsilon+\mathbf{S}(O,\delta)$ and
$\mathbf{B}_0+\mathbf{S}(O,\delta)$ into domains (part of conical
surface and part of spherical surface) and using \eqref{2.1} we get
the expressions for $\mu([\mathbf{B}_\varepsilon])$ and
$\mu([\mathbf{B}_0])$ in terms of flag density. The last expressions
we put into \eqref{2.1} and obtain \eqref{2.3}. The realization of this procedure can be found in \cite{Ara94} (see also \cite{Ara09}).

\section {AN INVERSION FORMULA FOR \eqref{1.5}}
\label{sec3}

To derive an inversion formula for \eqref{1.5}, we express $h$ in
terms of given function $\rho$.
Let $x\in e$, $<e,x>=\{f:\,\,f=(x,g,e)\}$ is the bundle of flags.

\begin{theorem}
\label{t1} Let $\rho$ be a $C^2$ smooth function on the space of
flags in $\mathbf R^3$, that is defined by transform \eqref{1.5}
of a sign measure $\mu$ in $\mathbf E$ with continuous density $h$.
For a given plane $e$ the following representation is valid
\begin{equation}\label{3.0}
h(e)=M(x)+\frac1{2\pi}\int_0^{2\pi}
\frac{\partial^2M(x)}{\partial^2_\varphi x}\,d\varphi-\frac2\pi
\int_0^{2\pi}(\rho''_{\phi\phi}(\varphi)+2\rho''_{\phi
y}(\varphi)+\rho''_{yy}(\varphi))d\varphi.
\end{equation}
where $x$ is a point on $e$, $\rho(\cdot)$ is the restriction of $\rho$ onto $<e,x>$, $M(x)$ is the bundle mass of $\rho$ at $x$.
\end{theorem}

\emph{Proof of the Theorem 6} Let $x\in e$ be a point. We consider a unit sphere
$\mathbf S^2_1$ tangent to $e$ at $x\in e$ and denote by
$A_n\subset\mathbf S^2_1$ the spherical disk centered at $x$ with
the spherical radius $\alpha_n$. By the mean-value theorem
\begin{equation}\label{3.1}
h(e)=\lim_{\alpha_n\to 0}\frac1{|A_n|}\int_{A_n}h^*(s)\,ds,
\end{equation}
where $|A_n|$ is the area of $A_n$.

\noindent On $\mathbf S^2_1$ we consider usual spherical coordinates
$(\varphi,\nu)$ with $x$ chosen for the pole. We consider the restriction of $M$ onto $\mathbf S^2_1$ and by $M(\varphi,\nu)$ we denote the value of the restriction at $(\varphi,\nu)\in\mathbf S^2_1$. Also, we consider the restriction of $\rho$ onto
$<\partial\mathrm{A}>$, where $A\subset \mathbf
S^2_1$ is the spherical disk centered at $x$ with the spherical
radius $\alpha$. By
$\rho(\varphi,\alpha)$ we denote the value of the restriction at $l=(\varphi,\alpha)\in\partial\mathrm{A}$.

\noindent We represent the integral in \eqref{3.1} according to
\eqref{2.3} and find the limit by decomposing the resulting terms in
powers of $\alpha_n$. We have
\begin{equation}\label{3.2}M(\varphi,\alpha_n)=
M(\varphi,0)+M'_\alpha(\varphi,0)\alpha_n+M''_{\alpha\alpha}(\varphi,0)
\frac{\alpha_n^2}2+o(\alpha_n^2).\end{equation} To decompose other
terms we need a lemma.

\begin{lemma}
\label{t5} For any differentiable flag function $\rho$
\begin{equation}\label{3.3}\frac{\partial\rho(\varphi,\alpha)}{\partial\alpha}
\mid_{\alpha=0}=\rho'_\phi(\varphi,0)+\rho'_y(\varphi,0).
\end{equation}
\end{lemma}

\noindent In its formulation $\rho'_y$ is defined assuming that the
positive normal of $e$ is parallel to the outer normal to $\mathbf
S^2_1$ at $x$. Note, that $\rho(\varphi) = \rho(\varphi,0)$ defined
on the bundle of flags $<e,x>$.

\noindent We have
\begin{equation}\label{3.4}\rho(\varphi,\alpha_n)=\rho(\varphi)+(\rho'_\phi(\varphi)+\rho'_y(\varphi)
) \alpha_n+(\rho''_{\phi\phi}(\varphi)+2\rho''_{\phi
y}(\varphi)+\rho''_{yy}(\varphi))\frac{\alpha_n^2}2+o(\alpha_n^2),\end{equation}
\begin{equation}\label{3.5}
\rho'_\phi(\varphi,\alpha_n)=\rho'_\phi(\varphi)+(\rho''_{\phi\phi}(\varphi)+\rho''_{\phi
y}(\varphi) ) \alpha_n+o(\alpha_n),\end{equation}
\begin{equation}\label{3.6}\rho'_y(\varphi,\alpha_n)=\rho'_y(\varphi)+(\rho''_{\phi y}(\varphi)+\rho''_{yy}(\varphi)
) \alpha_n+o(\alpha_n).\end{equation} After proper substitution and
taking into account that

\noindent $|A_n|=2\pi(1-\cos\alpha_n)\sim
\pi\alpha_n^2$, we find
\begin{multline}\label{3.7}
h(e)=M(x)+\frac{\partial M(x)}{\partial
n}+\\
+\frac1{2\pi}\int_0^{2\pi}M''_{\alpha\alpha}(\varphi,0)d\varphi-\frac2\pi
\int_{0}^{2\pi}(\rho''_{\phi\phi}(\varphi)+2\rho''_{\phi
y}(\varphi)+\rho''_{yy}(\varphi))d\varphi.\end{multline} Here
$\frac{\partial M(x)}{\partial n}$ is the derivative of $M$ at
$x\in e$, in the direction of normal to $e$ which is outer normal to
$\mathbf S^2_1$.

\noindent It is easy to prove that
\begin{equation}\label{3.8}M''_{\alpha\alpha}(\varphi,0)=
\frac{\partial^2M(x)}{\partial^2_\varphi x}- \frac{\partial
M(x)}{\partial n},\end{equation} where
$\frac{\partial^2M(x)}{\partial^2_\varphi x}$ denotes the second
derivative in the spatial direction which corresponds to the
direction $\varphi$ on $e$. We prefer to use this notation instead
of $M''_{\varphi\varphi}$, because the latter would correspond to
repeated differentiation on the spherical surface. After
substitution \eqref{3.8} into \eqref{3.7} we obtain \eqref{1.5}. Theorem 6 is proved.

The second integral in \eqref{1.5} contains the derivative
$\rho''_{\Phi\Phi}$. To present the
expression for $\rho'_\Phi$ in terms of the derivatives of
the function with respect to the parameters involved we introduce on $\mathbf S^2$ spherical coordinates in
which $\omega$ is specified by a pair $(\upsilon,\chi)$, where
$\upsilon$ denotes the latitude of $\omega$, and
$\chi$ denotes its longitude. The corresponding
derivatives with respect $\Phi$ are given by the expressions (see.
\cite{Ara09}):
\begin{equation}\label{3.9}
\varphi'_{\Phi}=-\tan\chi\cdot\sin\varphi\ \quad
\chi'_{\Phi}=-\cos\varphi\ ,\quad
\upsilon'_{\Phi}=\frac{\sin\varphi}{\cos\chi}\ .
\end{equation}
Using these formulas one can write the expression for
$\rho'_{\Phi}$.

\section {THE CONNECTION BETWEEN EQUATIONS \eqref{1.2} AND \eqref{1.5}}
\label{sec4}

For a given line $g$ and a point
$x\in g$, the corresponding bundle of flags we denote by
$$<x,g>=\{f:\,\,f=(x,g,e)\}.$$
Now we integrate \eqref{1.5} over all bundle of flags $<x,g>$ with respect to  Lebesgue measure $d\Phi$ on $\mathbf
S_\Omega$, where $\Omega$ is the direction of $g$ and $\mathbf
S_\Omega$ is the great circle on $\mathbf
S^2$ with pole at $\Omega$. We obtain
\begin{equation}\label{4.1}
\frac1{2\pi}\int_0^{2\pi}\rho(x,\Omega,\Phi)\,d\Phi=
\frac12\int_{\mathbf
S^2}|(\Omega,\xi)|\,h_{x}(\xi)\,d\xi,\end{equation}
since from the definition of flag density we have (see \cite{Amb90})
\begin{equation}\label{4.2}
\int_0^{2\pi}
\sin^2\alpha(\xi,f)\,d\,\Phi={2\pi}\,|(\Omega,\xi)|.
\end{equation}
Note that one can obtain \eqref{4.2} using the cosinus
theorem of spherical geometry.

It follows from \eqref{4.1} if $H$ is a smooth projective Finsler metric in $\mathbf R^3$ and $\rho$ is a flag density such that
\begin{equation}\label{4.3}
H(x,\Omega)=\frac1{\pi}
\int_{0}^{2\pi}\rho(x,\Omega,\Phi)\,d\Phi\,\,\,\,\emph{for}\,\,\,\,(x,\Omega)\in \mathbf R^3\times{\bold S}^2.\end{equation}
then equations \eqref{1.2} and \eqref{1.5} has the same unique solution.
Thus the problem of finding the solution of \eqref{1.2} we reduce to find flag density $\rho$ which satisfy \eqref{4.3} for a given smooth projective Finsler metric $H$.

Note that for a fixed $x\in\mathbf R^3$ the restriction $H(x,\cdot)$ of a smooth projective Finsler metric $H$ onto $\mathbf
S^2$ is uniquely determined by the restriction $h_{x}(\cdot)$ of the unique solution $h$ of \eqref{1.2} onto the bundle $[x]$. The restriction $h_{x}(\cdot)$ can be found by solution of the zonoid equation \eqref{1.4} for $H(\cdot)=H(x,\cdot)$.

Also, note that for a fixed $x\in\mathbf R^3$ the restriction $\rho(x,\cdot,\cdot)$
of the flag density $\rho$ of the sign measure $\mu(de)=h(e)de$ (where $h$ is the unique solution of \eqref{1.4}) is uniquely determined by the restriction $h_{x}(\cdot)$.

Hence for every fixed $x\in\mathbf R^3$, we
solve the zonoid equation for the function $H(x,\cdot)$. Since this function is of
class $C^\infty$, there exists a continuous even solution
$h(x,\cdot)$ on ${\bold S}^2$.
Then we put the function $h(x,\cdot)$ into equation \eqref{1.5} instead of $h_{x}(\cdot)$ and found the flag density $\rho(x,\cdot,\cdot)$ which depends on $x\in\mathbf R^3$. Thus we found the flag density $\rho$ which satisfy \eqref{4.3}.
According to Theorem 6, substituting the flag density $\rho$ into \eqref{3.1} we get the solution of \eqref{1.5} which coincides with the solution of \eqref{1.2}. Theorem 2 is proved.


\begin{thebibliography}{}

\bibitem{Amb76} R. V. Ambartzumian, A note on pseudo-metrics on the plane, {\em Z.Wahrsch. Verw. Gebiete},
vol. 37, pp. 145 – 155, 1976.


\bibitem{Amb87} R. V. Ambartzumian, Combinatorial integral geometry, metric and
zonoids, {\em Acta Appl. Math.}, vol.9, pp. 3 - 27, 1987.

\bibitem{Amb90} R. V. Ambartzumian, Factorization Calculus and Geometric
Probability,{\em Cambridge University Press}, 1990.

\bibitem{Amb94} R. V. Ambartzumian, V. K. Oganyan
Finite additive functionals in the space of planes, I, {\em Journal of
Contemporary Mathematical Analysis (Armenian Academy of Sciences)},
vol. 29, No. 4, pp. 3-51, 1994.

\bibitem{Amb98} R. V. Ambartzumian, V. K. Oganyan Parametric versions of Hilbert's fourth problem, {\em Israel J. Math.}, vol. 103, pp. 41–65, 1998.


\bibitem{Al} R. Alexander
Planes for which the lines are the shortest paths between points, {\em Illinois J. Math.}, vol 22, pp.177-190, 1978.

\bibitem{Alv03} J.C. Alvarez Paiva, Hilbert's fourth problem in two dimensions I, {\em Mass
Selecta: Teaching and Learning Advanced Undergraduate Mathematics} (S.
Katok, A. Sossinsky, S. Tabachnikov, eds.), pp. 165–183, Amer. Math. Soc.,
Providence, RI 2003.

\bibitem{Alv} J.C. Alvarez Paiva, Symplectic geometry and Hilbert's fourth problem, {\em J.
Differential Geometry} (to appear).

\bibitem{Alv98} J.C. Alvarez Paiva, E. Fernandes, Crofton formulas in projective Finsler
spaces, {\em Electron. Res. Announc. Amer. Math. Soc.}, vol. 4, pp. 91 – 100, 1998.


\bibitem{Ara92} R. H. Aramyan,  Recovering rose of directions
from flag densities in $\mathbf R^3$, {\em Journal of Contemporary
Mathematical Analysis (Armenian Academy of Sciences)}, vol. 27, No.
5, pp. 22 - 35, 1992.

\bibitem{Ara94} R. H. Aramyan, Measure generation in the space of planes and spherical Euler functionals, {\em Journal of Contemporary
Mathematical Analysis (Armenian Academy of Sciences)}, vol. 29, No.
4, pp. 52 - 74, 1994.

\bibitem{Ara09} R. Aramyan,
Solution of one integral equation on the sphere by methods of
integral geometry, {\em Doklady Mathematics,} vol. 79, No. 3, pp.
325-328, 2009 [{\em Doklady Akademii Nauk (Russian Academy of
Sciences),} vol. 426, No. 2, pp. 151-155, 2009].

\bibitem{Ara10} R. Aramyan, Measure of planes intersecting a convex body, {\em Sutra: Inter. J. of Math. Science}, vol. 3 (1), pp. 1 - 7, 2010.

\bibitem{SW} R. Schneider, W.Weil, Zonoids and related topics, {\em Convexity and Its Applications},
(P.M. Gruber, J.M. Wills, eds.), pp. 296–317, Birkh.auser, Basel
1983.

\bibitem{Bu} H. Busemann, Problem IV: Desarguesian spaces in Mathematical developments arising from Hilbert problems, {\em Proc. Sympos. Pure Math.}, Vol. 28, Amer. Math.
Soc., Providence, R. I., 1976.

\bibitem{Pog} A. V. Pogorelov,  Hilbert's Fourth Problem, {\em Wilson and Sons}, 1979 (Russian Original: Izdat. "Nauka", Moscow
1974).

\bibitem{Pa}G. Yu. Panina, Convex bodies and translation invariant measures,
{\em Zap. Nauchn. Semin. LOMI}, vol. 157, pp. 1986.

\bibitem{Sch67} R. Schneider, Zu einem Problem von Shephard \"uber die Projektionen konvexer
K\"orper, {\em Math. Z.}, vol. 191, pp. 71 – 82, 1967.

\bibitem{Sch02}  R. Schneider, On integral geometry in projective Finsler spaces, {\em Journal of Contemporary
Mathematical Analysis (Armenian Academy of Sciences)}, vol. 37, pp. 34–51, 2002.

\bibitem{Sch05}  R. Schneider, Crofton measures in projective Finsler spaces,
{\em In: Integral Geometry and Convexity (Proc. Int. Conf., Wuhan, China, Oct. 2004; eds. E.L. Grinberg, S. Li, G. Zhang, J. Zhou)}, World Scientific, New Jersey,  pp. 67 - 98, 2006.


\end{thebibliography}
\end{document}